\documentclass[12pt]{amsart}

\usepackage[T1]{fontenc}
\usepackage{lmodern}
\usepackage{amsmath,amssymb,amsthm,mathtools}
\usepackage[a4paper,margin=1.05in]{geometry}
\usepackage{microtype}
\usepackage{needspace}
\usepackage{cite}
\usepackage{xcolor}
\usepackage{hyperref}
\hypersetup{colorlinks=true,linkcolor=blue,citecolor=blue,urlcolor=blue,
  pdftitle={Bergman Kernel Zeros on Rational Hartogs Triangles},
  pdfauthor={Qian Fu}}
\newtheorem{theorem}{Theorem}[section]
\newtheorem*{emconjecture}{Edholm--Mathew conjecture}
\newtheorem{proposition}[theorem]{Proposition}
\newtheorem{lemma}[theorem]{Lemma}
\newtheorem{corollary}[theorem]{Corollary}
\theoremstyle{definition}
\newtheorem{remark}[theorem]{Remark}
\numberwithin{equation}{section}

\makeatletter
\renewenvironment{proof}[1][\proofname]{\par
  \pushQED{\qed}%
  \normalfont \topsep6\p@\@plus6\p@\relax
  \trivlist
  \item[\hskip\labelsep\bfseries #1\@addpunct{.}]\ignorespaces
}{%
  \popQED\endtrivlist\@endpefalse
}
\makeatother

\newcommand{\C}{\mathbb C}
\newcommand{\D}{\mathbb D}
\newcommand{\Dstar}{\mathbb D^{*}}
\newcommand{\R}{\mathbb R}
\newcommand{\Q}{\mathbb Q}
\newcommand{\Z}{\mathbb Z}

\title[Bergman kernel zeros on rational Hartogs triangles]
{Bergman Kernel Zeros on Rational Hartogs Triangles}

\author[G. T. Deng]{Guantie Deng}
\address{School of Mathematical Sciences, Beijing Normal University, Beijing 100875, P.R. China.}\email{96022@bnu.edu.cn}

\author[Q. Fu]{Qian Fu$^{*}$}
\address{Teaching Experiment Platform, Beijing Normal University, Zhuhai 519087, P. R. China}
\email{qianfu@bnu.edu.cn}

\date{}
\thanks{${ }^{*}$Corresponding author.}
\begin{document}

\begin{abstract}
We prove the Edholm--Mathew nonvanishing conjecture for the palindromic
polynomials associated with rational Hartogs triangles.
For relatively prime positive integers $m>n$, we show that
$e^{-i(m-n)\theta}Q_{m,n}(e^{i\theta})$ is strictly positive for every
real $\theta$ and derive an explicit lower bound. The proof expresses
this trigonometric polynomial as a finite average of products of
unnormalized Fej\'er kernels. It follows that $Q_{m,n}$ has exactly
$m-n$ zeros in the punctured unit disc, counted with multiplicity.
As a consequence, a rational Hartogs triangle has a zero-free Bergman
kernel if and only if its exponent is the reciprocal of a positive
integer.
\end{abstract}

\maketitle

\section{Introduction}\label{sec:intro}

Let $\Omega$ be a bounded domain in $\C^d$, and let
$A^2(\Omega)=\mathcal O(\Omega)\cap L^2(\Omega)$, with respect to
Lebesgue volume. The orthogonal projection onto $A^2(\Omega)$ has the
integral representation
\[
 P_\Omega f(z)=\int_\Omega K_\Omega(z,w)f(w)\,dV(w),
 \qquad f\in L^2(\Omega).
\]
Although $K_\Omega(z,z)>0$ for every $z\in\Omega$, the kernel may
vanish away from the diagonal. The Lu Qi-Keng problem asks which
domains have a zero-free Bergman kernel on $\Omega\times\Omega$;
such domains are called \emph{Lu Qi-Keng domains}. This question
arises in the study of Bergman representative coordinates and is
invariant under biholomorphic mappings; see Boas \cite{Boas2000}.
Explicit examples show that even convexity does not prevent kernel
zeros \cite{BoasFuStraube1999}. Thus, determining the zero set requires
information beyond diagonal positivity or geometric convexity.

We consider the generalized Hartogs triangles
\[
 H_\gamma=\{(z_1,z_2)\in\C^2:|z_1|^\gamma<|z_2|<1\},
 \qquad \gamma>0,
\]
and write $K_\gamma=K_{H_\gamma}$. For a rational exponent in lowest
terms,
\begin{equation}\label{eq:domain}
 H_{m/n}=\{(z_1,z_2)\in\C^2:|z_1|^m<|z_2|^n<1\},
 \qquad m,n\in\Z_+,\quad \gcd(m,n)=1,
\end{equation}
where $\Z_+$ denotes the positive integers. These bounded pseudoconvex
Reinhardt domains have singular boundary and exhibit a close connection
between the arithmetic of the exponent and the analytic properties of
the Bergman projection. Chakrabarti and Zeytuncu
\cite{ChakrabartiZeytuncu2016} established sharp weighted estimates on
the classical Hartogs triangle. Edholm and McNeal
\cite{EdholmMcNeal2017} proved that the Bergman projection on $H_{m/n}$
is bounded on $L^p$ precisely when
\[
 \frac{2(m+n)}{m+n+1}<p<\frac{2(m+n)}{m+n-1},
\]
whereas for irrational $\gamma$ it is bounded only for $p=2$.

Subsequent work has developed several aspects of this operator theory.
Huo and Wick obtained weighted norm estimates in terms of
Bekoll\'e--Bonami-type constants \cite{HuoWickJFA2020} and determined
the weak-type behavior at the endpoints of the boundedness interval on
the classical Hartogs triangle \cite{HuoWickBLMS2020}.
Chakrabarti, Edholm, and McNeal \cite{ChakrabartiEdholmMcNeal2019}
studied the consequences of irregularity for duality and approximation
in Bergman spaces. More recently, Chakrabarti and Edholm
\cite{ChakrabartiEdholm2024} constructed monomial basis projections
adapted to $L^p$-Bergman spaces on Reinhardt domains. These results
emphasize the role of the monomial structure of the kernel. The present
paper uses that structure to study its zeros.

For rational exponents, a complementary description is provided by
finite holomorphic quotients. Write
$\D=\{\zeta\in\C:|\zeta|<1\}$ and $\Dstar=\D\setminus\{0\}$.
The monomial map
\[
 \pi_{m,n}:\D\times\Dstar\longrightarrow H_{m/n},
 \qquad \pi_{m,n}(u,v)=(uv^n,v^m),
\]
is an $m$-sheeted proper holomorphic covering. Its cyclic deck group
acts by
\[
 (u,v)\longmapsto(\omega^{-n}u,\omega v),\qquad \omega^m=1.
\]
Bell's transformation formula \cite{Bell1982} relates the kernel on the
quotient to the kernel on the product domain. This approach underlies
the $L^p$ theory of quotient domains developed by Bender, Chakrabarti,
Edholm, and Mainkar \cite{BenderChakrabartiEdholmMainkar2022}, as well as
explicit kernel formulas for elementary Reinhardt domains
\cite{ChakrabartiKonkelMainkarMiller2020}, two-dimensional monomial
polyhedra \cite{Almughrabi2023}, and general monomial polyhedra
\cite{ChakrabartiCinzoriGaidhaneGregoryWright2024}. Although
$\D\times\Dstar$ has a zero-free kernel, zeros can appear on its
finite quotients.

For the Lu Qi-Keng problem on Hartogs triangles, Edholm
\cite{Edholm2016} proved that $H_j$ has kernel zeros for each integer
$j\ge2$, while $H_{1/j}$ is a Lu Qi-Keng domain for every
$j\in\Z_+$. Edholm and McNeal \cite{EdholmMcNeal2017} extended the
existence of zeros to every real $\gamma>2$. For $\gamma>1$, their
calculation on the coordinate slice $z_1\overline{w_1}=0$ gives
\[
 K_\gamma\bigl((0,z_2),(0,w_2)\bigr)
 =\frac{1+(\gamma-1)t}{\gamma\pi^2t(1-t)^2},
 \qquad t=z_2\overline{w_2}.
\]
The numerator vanishes at $t=-1/(\gamma-1)$, which belongs to
$\Dstar$ exactly when $\gamma>2$. The interval $1<\gamma<2$
therefore requires a different choice of slice.

Edholm and Mathew \cite{EdholmMathew2026} introduced the oblique slice
$s=t$, where
\[
 s=z_1\overline{w_1},\qquad t=z_2\overline{w_2}.
\]
For relatively prime $m>n$, the formula of Almughrabi
\cite[Proposition~2]{Almughrabi2023}, also recalled in
\cite[Proposition~1.10]{EdholmMathew2026}, is
\begin{equation}\label{eq:kernel}
 K_{m/n}(z,w)
 =\frac{P_{m,n}(s,t)}{m\pi^2(1-t)^2(t^n-s^m)^2},
\end{equation}
where $P_{m,n}$ is given in \eqref{eq:numerator} below. Following
Edholm and Mathew, set
\begin{equation}\label{eq:Q}
 k=m-n,\qquad Q_{m,n}(\zeta)=\zeta^{1-2n}P_{m,n}(\zeta,\zeta).
\end{equation}
They proved that $Q_{m,n}$ is a polynomial of degree $2k$ with positive
coefficients and palindromic symmetry
\cite[Theorem~1.5 and Proposition~3.22]{EdholmMathew2026}:
\begin{equation}\label{eq:reciprocal}
 Q_{m,n}(\zeta)=\zeta^{2k}Q_{m,n}(\zeta^{-1}),\qquad \zeta\ne0.
\end{equation}
Since $(\zeta,\zeta)\in H_{m/n}$ for $0<|\zeta|<1$, roots of
$Q_{m,n}$ in $\Dstar$ produce zeros of $K_{m/n}$. The reciprocal
symmetry pairs roots inside and outside the unit circle, but does not
by itself exclude roots on the circle.

Edholm and Mathew proved that $Q_{m,n}$ has no unit-circle roots when
$k=1$ or $2$ \cite[Theorem~1.7]{EdholmMathew2026}. In particular, they
obtained rational exponents decreasing to $1$ for which the kernels
have zeros, despite the nonvanishing of $K_1$. In
\cite[Section~4.3]{EdholmMathew2026}, they also announced the cases
$k=3,4,6$ and infinitely many pairs for each fixed $k$.
Their general question is the following
\cite[Conjecture~4.20]{EdholmMathew2026}.

\begin{emconjecture}
Let $m>n$ be relatively prime positive integers. Then the polynomial
$Q_{m,n}$ has no zeros on the unit circle. Consequently, $Q_{m,n}$
has exactly $m-n$ zeros in $\Dstar$, counted with multiplicity.
\end{emconjecture}

We prove this conjecture by establishing a stronger positivity
property on the unit circle.

\begin{theorem}\label{thm:positive}
Let $m>n$ be relatively prime positive integers and set $k=m-n$. Then
\begin{equation}\label{eq:strict-positive}
 e^{-ik\theta}Q_{m,n}(e^{i\theta})>0
 \qquad\text{for every }\theta\in\R.
\end{equation}
More precisely,
\begin{equation}\label{eq:lower-bound}
 e^{-ik\theta}Q_{m,n}(e^{i\theta})
 \ge m\sin^2(\theta/2),
 \qquad Q_{m,n}(1)=m^3.
\end{equation}
\end{theorem}

Theorem~\ref{thm:positive} determines the number of zeros of
$Q_{m,n}$ in $\Dstar$ for every relatively prime pair $m>n$.
Combining this with the monomial biholomorphisms of generalized Hartogs
triangles \cite[Theorem~2.19]{EdholmMathew2026} gives a complete
classification for rational exponents.

\begin{corollary}\label{cor:rational-classification}
Let $\gamma=m/n>0$, where $m,n\in\Z_+$ and $\gcd(m,n)=1$.
The domain $H_\gamma$ is a Lu Qi-Keng domain if and only if $m=1$.
Equivalently, for $\gamma\in\Q_{>0}$,
\begin{equation}\label{eq:classification}
 K_\gamma\text{ is zero-free on }H_\gamma\times H_\gamma
 \quad\Longleftrightarrow\quad
 \gamma=\frac1j\text{ for some }j\in\Z_+.
\end{equation}
\end{corollary}

The main step is a positive representation of
$e^{-ik\theta}Q_{m,n}(e^{i\theta})$. Rather than expanding the
coefficients separately for each value of $k$, we retain the product
of triangular weights in Almughrabi's formula. After centering, each
weight is the discrete autocorrelation of the sequence consisting of
$m$ ones. Its Fourier polynomial is therefore the squared modulus of
a finite geometric sum, namely the unnormalized Fej\'er kernel.
The two weights are coupled by one congruence modulo $m$.
Finite Fourier orthogonality separates this congruence and gives the
product representation in Lemma~\ref{lem:filter}. Coprimality yields
an exact identity for the sum of the first factors, while the second
factors admit a common positive lower bound away from $2\pi\Z$.
The remaining angles are evaluated directly. This proves
Theorem~\ref{thm:positive} for all $m>n$ without a case distinction
in $m-n$.

Section~\ref{sec:coefficients} establishes the centered coefficient
formula. Section~\ref{sec:fejer} proves
Theorem~\ref{thm:positive} and the Edholm--Mathew conjecture, and
records the equivalent positivity statement on a real interval.
Section~\ref{sec:zeros} constructs kernel zeros and proves
Corollary~\ref{cor:rational-classification}.

\section{The kernel formula and centered coefficients}\label{sec:coefficients}

Throughout this section and Section~\ref{sec:fejer}, let $m>n$ be
relatively prime positive integers and put $k=m-n$. For $v\in\Z$, define
\begin{equation}\label{eq:D}
 D_m(v)=
 \begin{cases}
  v+1,&0\le v\le m-1,\\
  2m-1-v,&m\le v\le2m-2,\\
  0,&\text{otherwise}.
 \end{cases}
\end{equation}
Recall the kernel formula
\[
 K_{m/n}(z,w)
 =\frac{P_{m,n}(s,t)}{m\pi^2(1-t)^2(t^n-s^m)^2},
 \qquad s=z_1\overline{w_1},\quad t=z_2\overline{w_2}.
\]
Its numerator is
\begin{equation}\label{eq:numerator}
 P_{m,n}(s,t)
 =\sum_{a=0}^{2m-2}\sum_{b=0}^{2n}
 D_m(a)D_m(mb+na+m+n-1-2mn)s^at^b;
\end{equation}
see \cite[Proposition~2]{Almughrabi2023} and
\cite[Proposition~1.10]{EdholmMathew2026}. 

For $z,w\in H_{m/n}$, their coordinate products satisfy
\begin{equation}\label{eq:products}
 |s|^m<|t|^n<1.
\end{equation}
Thus neither $1-t$ nor $t^n-s^m$ vanishes in \eqref{eq:kernel}.
Conversely, every $(s,t)\in H_{m/n}$ is obtained as such a pair of
coordinate products. Indeed, writing
$s=r_1e^{i\phi_1}$ and $t=r_2e^{i\phi_2}$, take
\[
 z=(\sqrt{r_1},\sqrt{r_2}),\qquad
 w=(\sqrt{r_1}e^{-i\phi_1},\sqrt{r_2}e^{-i\phi_2}).
\]
When $s=0$, the choice of $\phi_1$ is immaterial. The defining inequality
for both points is the square root of \eqref{eq:products}. This is the
special case for Hartogs triangles of the factorization discussed in
\cite[Section~2.2.2]{EdholmMathew2026}.

We center the coefficients in \eqref{eq:D} by setting
\begin{equation}\label{eq:W}
 W_m(A)=
 \begin{cases}
  m-|A|,&|A|\le m-1,\\
  0,&|A|\ge m,
 \end{cases}
 \qquad A\in\Z.
\end{equation}
Then $D_m(v)=W_m(v-(m-1))$.

\begin{lemma}\label{lem:center}
For $\zeta\ne0$, let
$Q_{m,n}(\zeta)=\zeta^{1-2n}P_{m,n}(\zeta,\zeta)$ as in
\eqref{eq:Q}. Then
\begin{equation}\label{eq:center}
 Q_{m,n}(\zeta)
 =\sum_{A,B\in\Z}W_m(A)W_m(nA+mB)\zeta^{k+A+B}.
\end{equation}
The sum is finite and extends to a polynomial of degree $2k$ with real
nonnegative coefficients and positive constant and leading coefficients.
It also satisfies the reciprocal identity
\[
 Q_{m,n}(\zeta)=\zeta^{2k}Q_{m,n}(\zeta^{-1}),\qquad \zeta\ne0.
\]
\end{lemma}

\begin{proof}
Set
\[
 A=a-(m-1),\qquad B=b-n,\qquad V=nA+mB.
\]
Using $D_m(v)=W_m(v-(m-1))$, we obtain
\begin{align*}
 D_m(a)&=W_m(A),\\
 D_m(mb+na+m+n-1-2mn)
 &=W_m\bigl(mb+na+n-2mn\bigr)\\
 &=W_m\bigl(m(b-n)+n(a-m+1)\bigr)
 =W_m(V),\\
 1-2n+a+b&=k+A+B.
\end{align*}
Moreover,
\[
 \begin{aligned}
 0\le a\le2m-2&\iff |A|\le m-1,\\
 0\le b\le2n&\iff |B|\le n.
 \end{aligned}
\]
Consequently, for $\zeta\ne0$,
\begin{equation}\label{eq:center-finite}
 Q_{m,n}(\zeta)
 =\sum_{\substack{|A|\le m-1\\|B|\le n}}
 W_m(A)W_m(V)\zeta^{k+A+B}.
\end{equation}
For any nonzero summand,
\begin{align*}
 W_m(A)W_m(V)\ne0
 &\implies |A|\le m-1,\quad |V|\le m-1,\\
 &\implies m|B|=|V-nA|\\
 &\hspace{3.7em}\le |V|+n|A|
 \le(n+1)(m-1)<m(n+1),\\
 &\implies |B|\le n\qquad(B\in\Z).
\end{align*}
Thus the sum in \eqref{eq:center-finite} may be extended to
$\Z^2$, giving \eqref{eq:center}; it remains finite since
\[
 W_m(A)W_m(nA+mB)\ne0
 \implies |A|\le m-1,\quad |B|\le n.
\]
Next, the identity
\begin{equation}\label{eq:linear-relation}
 m(A+B)=kA+V
\end{equation}
gives, for every nonzero summand,
\begin{align*}
 m|A+B|&=|kA+V|\\
 &\le k|A|+|V|\\
 &\le(k+1)(m-1)<m(k+1),\\
 -k\le A+B&\le k\qquad(A+B\in\Z),\\
 0\le k+A+B&\le2k.
\end{align*}
Hence $Q_{m,n}\in\R[\zeta]$, all its coefficients are
nonnegative, and $\deg Q_{m,n}\le2k$.

\Needspace{9\baselineskip}
For
\[
 A_0=m-1,\qquad B_0=1-n,
\]
we have
\begin{align*}
 A_0+B_0&=k,\\
 nA_0+mB_0&=n(m-1)+m(1-n)=k,\\
 W_m(A_0)W_m(nA_0+mB_0)
 &=W_m(m-1)W_m(k)=m-k=n>0.
\end{align*}
Thus, if $Q_{m,n}(\zeta)=\sum_{j=0}^{2k}c_j\zeta^j$, then
\[
 c_{2k}\ge n>0,
 \qquad \deg Q_{m,n}=2k.
\]
Finally, since $W_m(-x)=W_m(x)$,
\begin{align*}
 \zeta^{2k}Q_{m,n}(\zeta^{-1})
 &=\sum_{A,B\in\Z}W_m(A)W_m(nA+mB)
   \zeta^{k-A-B}\\
 &=\sum_{A,B\in\Z}W_m(-A)W_m(-nA-mB)
   \zeta^{k+A+B}\\
 &=Q_{m,n}(\zeta).
\end{align*}
Therefore
\[
 c_j=c_{2k-j}\quad(0\le j\le2k),
 \qquad c_0=c_{2k}\ge n>0.
\]
\end{proof}

\section{Proof of Theorem~\ref{thm:positive}}\label{sec:fejer}

For $u\in\R$, let
\begin{equation}\label{eq:fejer}
 T_m(u)=\sum_{A=-m+1}^{m-1}W_m(A)e^{iAu}
       =\left|\sum_{j=0}^{m-1}e^{iju}\right|^2.
\end{equation}
Indeed, for $A\in\Z$ and $0\le j,r\le m-1$, the condition
$j-r=A$ gives
\[
\begin{aligned}
 &\#\{(j,r)\in\{0,\ldots,m-1\}^2:j-r=A\}\\
 &\quad=
 \begin{cases}
  \#\{r:0\le r\le m-1-A\}=m-A,&0\le A\le m-1,\\
  \#\{r:-A\le r\le m-1\}=m+A,&-m+1\le A<0,\\
  0,&|A|\ge m
 \end{cases}\\
 &\quad=W_m(A).
\end{aligned}
\]
Consequently,
\[
\begin{aligned}
 \left|\sum_{j=0}^{m-1}e^{iju}\right|^2
 &=\sum_{j=0}^{m-1}\sum_{r=0}^{m-1}e^{i(j-r)u}\\
 &=\sum_{A=-(m-1)}^{m-1}
   \#\{(j,r):0\le j,r\le m-1,\ j-r=A\}\,e^{iAu}\\
 &=\sum_{A=-(m-1)}^{m-1}W_m(A)e^{iAu}=T_m(u).
\end{aligned}
\]
Thus $T_m(u)\ge0$. More explicitly,
\begin{equation}\label{eq:fejer-sine}
 T_m(u)=
 \begin{cases}
  \displaystyle\frac{\sin^2(mu/2)}{\sin^2(u/2)},&u\notin2\pi\Z,\\[6pt]
  m^2,&u\in2\pi\Z.
 \end{cases}
\end{equation}
This is the unnormalized Fej\'er kernel; see, for example,
\cite[Chapter~III]{Zygmund2002}.

Define
\begin{equation}\label{eq:F}
 F_{m,n}(\theta)
 =\sum_{A,B\in\Z}W_m(A)W_m(nA+mB)e^{i(A+B)\theta}.
\end{equation}
By Lemma~\ref{lem:center},
\begin{equation}\label{eq:F-Q}
 F_{m,n}(\theta)=e^{-ik\theta}Q_{m,n}(e^{i\theta}).
\end{equation}
In particular, $F_{m,n}$ is a $2\pi$-periodic trigonometric polynomial.

\Needspace{11\baselineskip}
\begin{lemma}\label{lem:filter}
For every $\theta\in\R$,
\begin{equation}\label{eq:filter}
 F_{m,n}(\theta)
 =\frac1m\sum_{\ell=0}^{m-1}
 T_m\!\left(\frac{k\theta-2\pi n\ell}{m}\right)
 T_m\!\left(\frac{\theta+2\pi\ell}{m}\right).
\end{equation}
Moreover,
\begin{equation}\label{eq:sum-first}
 \sum_{\ell=0}^{m-1}
 T_m\!\left(\frac{k\theta-2\pi n\ell}{m}\right)=m^2.
\end{equation}
\end{lemma}

\begin{proof}
Put
\[
 V=nA+mB,\qquad B=\frac{V-nA}{m}.
\]
The change of variables $(A,B)\mapsto(A,V)$ is a bijection from
$\Z^2$ onto
\[
 \{(A,V)\in\Z^2:V\equiv nA\pmod m\}.
\]
Since $k=m-n$,
\[
 A+B=A+\frac{V-nA}{m}=\frac{kA+V}{m},
\]
and hence \eqref{eq:F} becomes
\begin{equation}\label{eq:F-congruence}
 F_{m,n}(\theta)
 =\sum_{\substack{A,V\in\Z\\V\equiv nA\ ({\rm mod}\ m)}}
 W_m(A)W_m(V)e^{i(kA+V)\theta/m}.
\end{equation}
For every integer $q$, the geometric-sum formula gives
\begin{equation}\label{eq:roots-filter}
 \sum_{\ell=0}^{m-1}e^{2\pi i\ell q/m}
 =\begin{cases}
 m, &m\mid q,\\[3pt]
 \displaystyle\frac{1-e^{2\pi iq}}{1-e^{2\pi iq/m}}=0,
    &m\nmid q.
 \end{cases}
\end{equation}
Indeed, if $m\mid q$, each summand equals $1$; otherwise,
$e^{2\pi iq/m}\ne1$ and $e^{2\pi iq}=1$.
Taking $q=V-nA$ in \eqref{eq:roots-filter}, we obtain
\begin{align*}
 F_{m,n}(\theta)
 &=\frac1m\sum_{A,V\in\Z}
 W_m(A)W_m(V)e^{i(kA+V)\theta/m}
 \sum_{\ell=0}^{m-1}e^{2\pi i\ell(V-nA)/m}\\
 &=\frac1m\sum_{\ell=0}^{m-1}
 \sum_{A,V\in\Z}W_m(A)W_m(V)
 e^{iA(k\theta-2\pi n\ell)/m}
 e^{iV(\theta+2\pi\ell)/m}\\
 &=\frac1m\sum_{\ell=0}^{m-1}
 \left(\sum_{A\in\Z}W_m(A)
 e^{iA(k\theta-2\pi n\ell)/m}\right)\\
 &\hspace{6em}\cdot
 \left(\sum_{V\in\Z}W_m(V)
 e^{iV(\theta+2\pi\ell)/m}\right)\\
 &=\frac1m\sum_{\ell=0}^{m-1}
 T_m\!\left(\frac{k\theta-2\pi n\ell}{m}\right)
 T_m\!\left(\frac{\theta+2\pi\ell}{m}\right).
\end{align*}
All rearrangements are valid because $W_m$ has finite support.
This proves \eqref{eq:filter}.

To establish \eqref{eq:sum-first}, expand the first Fej\'er factor:
\begin{align*}
 &\sum_{\ell=0}^{m-1}
 T_m\!\left(\frac{k\theta-2\pi n\ell}{m}\right)\\
 &\quad=\sum_{\ell=0}^{m-1}
 \sum_{A=-m+1}^{m-1}W_m(A)
 e^{iA(k\theta-2\pi n\ell)/m}\\
 &\quad=\sum_{A=-m+1}^{m-1}W_m(A)e^{ikA\theta/m}
 \sum_{\ell=0}^{m-1}e^{-2\pi inA\ell/m}\\
 &\quad=m\sum_{\substack{-m+1\le A\le m-1\\m\mid nA}}
 W_m(A)e^{ikA\theta/m}\\
 &\quad=mW_m(0)=m^2.
\end{align*}
In the last step, \eqref{eq:roots-filter} and $\gcd(m,n)=1$
were used together with
\[
 \bigl(m\mid nA\ \text{and}\ |A|\le m-1\bigr)
 \iff
 \bigl(m\mid A\ \text{and}\ |A|\le m-1\bigr)
 \iff A=0.
\]
\end{proof}

\begin{proposition}\label{prop:strict}
The trigonometric polynomial $F_{m,n}$ is strictly positive on $\R$.
It satisfies
\begin{equation}\label{eq:F-bound}
 F_{m,n}(\theta)\ge m\sin^2(\theta/2)
 \qquad(\theta\in\R),
\end{equation}
and $F_{m,n}(2\pi j)=m^3$ for every $j\in\Z$.
\end{proposition}

\begin{proof}
For $0\le\ell\le m-1$, put
\[
 u_\ell=\frac{k\theta-2\pi n\ell}{m},\qquad
 v_\ell=\frac{\theta+2\pi\ell}{m}.
\]
Recall from \eqref{eq:fejer} that
\[
 T_m(u)=\left|\sum_{j=0}^{m-1}e^{iju}\right|^2\ge0,
 \qquad u\in\R.
\]
In the present notation, the two identities in Lemma~\ref{lem:filter}
read
\begin{equation}\label{eq:positive-factors}
 \begin{aligned}
  \sum_{\ell=0}^{m-1}T_m(u_\ell)&=m^2,\\
  F_{m,n}(\theta)&=\frac1m\sum_{\ell=0}^{m-1}
                T_m(u_\ell)T_m(v_\ell).
 \end{aligned}
\end{equation}

Suppose first that $\theta\notin2\pi\Z$. Then, for every
$0\le\ell\le m-1$,
\[
 e^{imv_\ell}=e^{i\theta}\ne1,
 \qquad e^{iv_\ell}\ne1,
 \qquad \sin^2(v_\ell/2)>0.
\]
The geometric-sum identity and \eqref{eq:fejer-sine} give
\begin{align}
 T_m(v_\ell)
 &=\left|\sum_{j=0}^{m-1}e^{ijv_\ell}\right|^2
 =\left|\frac{1-e^{imv_\ell}}{1-e^{iv_\ell}}\right|^2
 \notag\\
 &=\frac{\sin^2((\theta+2\pi\ell)/2)}
 {\sin^2((\theta+2\pi\ell)/(2m))}
 =\frac{\sin^2(\theta/2)}{\sin^2(v_\ell/2)}
 \ge\sin^2(\theta/2)>0.
 \label{eq:second-factor}
\end{align}
Here $0<\sin^2(v_\ell/2)\le1$. Consequently,
\begin{align*}
 F_{m,n}(\theta)
 &=\frac1m\sum_{\ell=0}^{m-1}T_m(u_\ell)T_m(v_\ell)\\
 &\ge\frac{\sin^2(\theta/2)}m\sum_{\ell=0}^{m-1}T_m(u_\ell)\\
 &=m\sin^2(\theta/2)>0.
\end{align*}

Now let $\theta\in2\pi\Z$. Every frequency $A+B$ in the
finite sum defining $F_{m,n}$ is an integer, and hence
\[
 e^{i(A+B)(\theta+2\pi j)}=e^{i(A+B)\theta},\qquad
 F_{m,n}(\theta+2\pi j)=F_{m,n}(\theta)\quad(j\in\Z).
\]
It therefore suffices to take $\theta=0$.
For $1\le\ell\le m-1$, the geometric-sum formula yields
\[
 T_m\!\left(\frac{2\pi\ell}{m}\right)
 =\left|\sum_{j=0}^{m-1}e^{2\pi ij\ell/m}\right|^2
 =\left|\frac{1-e^{2\pi i\ell}}
 {1-e^{2\pi i\ell/m}}\right|^2=0,
\]
whereas $T_m(0)=m^2$. In \eqref{eq:filter} at $\theta=0$,
only the summand $\ell=0$ remains, and hence
\begin{align*}
 F_{m,n}(0)
 &=\frac1m\sum_{\ell=0}^{m-1}
 T_m\!\left(-\frac{2\pi n\ell}{m}\right)
 T_m\!\left(\frac{2\pi\ell}{m}\right)\\
 &=\frac1m T_m(0)^2=\frac1m(m^2)^2=m^3>0.
\end{align*}
Thus $F_{m,n}(2\pi j)=m^3$ for all $j\in\Z$, while
$m\sin^2(\pi j)=0$. Combining the two cases proves both strict
positivity and \eqref{eq:F-bound} on $\R$.
\end{proof}

\Needspace{11\baselineskip}
\begin{proof}[Proof of Theorem~\ref{thm:positive}]
Set $k=m-n$. Evaluating the centered expression for $Q_{m,n}$
from Lemma~\ref{lem:center} at $\zeta=e^{i\theta}$, we obtain
\[
\begin{aligned}
 e^{-ik\theta}Q_{m,n}(e^{i\theta})
 &=\sum_{A,B\in\Z}W_m(A)W_m(nA+mB)e^{i(A+B)\theta}\\
 &=F_{m,n}(\theta).
\end{aligned}
\]
Thus the phase-adjusted boundary value of $Q_{m,n}$ is exactly the
real trigonometric polynomial studied in Proposition~\ref{prop:strict}.
The proposition establishes its strict positivity and the lower bound
\[
 e^{-ik\theta}Q_{m,n}(e^{i\theta})
 =F_{m,n}(\theta)\ge m\sin^2(\theta/2),
 \qquad F_{m,n}(\theta)>0\quad(\theta\in\R).
\]
Finally, taking $\theta=0$ and using $F_{m,n}(0)=m^3$, we get
\[
 Q_{m,n}(1)=F_{m,n}(0)=m^3,
\]
which proves all the assertions.
\end{proof}

\begin{proof}[Proof of the Edholm--Mathew conjecture]
Put $k=m-n$. Theorem~\ref{thm:positive} gives
\[
 Q_{m,n}(e^{i\theta})=e^{ik\theta}F_{m,n}(\theta)\ne0
 \qquad(\theta\in\R),
\]
since $F_{m,n}(\theta)>0$. Thus there are no roots on the unit circle.
Recall from Lemma~\ref{lem:center} that
\[
 \deg Q_{m,n}=2k,\qquad Q_{m,n}(0)>0,\qquad
 Q_{m,n}(\zeta)=\zeta^{2k}Q_{m,n}(\zeta^{-1})\quad(\zeta\ne0).
\]
In particular, every root is nonzero.

Let $a$ be a root of multiplicity $d$. Write
\[
 Q_{m,n}(\zeta)=(\zeta-a)^dG(\zeta),\qquad G(a)\ne0.
\]
The reciprocal identity gives, for $\zeta\ne0$,
\begin{align*}
 Q_{m,n}(\zeta)
 &=\zeta^{2k}(\zeta^{-1}-a)^dG(\zeta^{-1})\\
 &=(-a)^d\zeta^{2k-d}(\zeta-a^{-1})^dG(\zeta^{-1}).
\end{align*}
The factor multiplying $(\zeta-a^{-1})^d$ is holomorphic near
$a^{-1}$, and its value there is
\[
 (-a)^d(a^{-1})^{2k-d}G(a)\ne0.
\]
Consequently, $a^{-1}$ is a root of the same multiplicity $d$.
The map $a\mapsto a^{-1}$ is an involution, and
\[
 0<|a|<1\quad\Longleftrightarrow\quad |a^{-1}|>1.
\]
It therefore gives a multiplicity-preserving bijection between the
roots inside and outside the unit circle. If $N_{<}$ and $N_{>}$
denote their respective numbers, counted with multiplicity, the
fundamental theorem of algebra gives
\[
 N_{<}=N_{>},\qquad N_{<}+N_{>}=2k,
 \qquad N_{<}=N_{>}=k=m-n.
\]
Since $Q_{m,n}(0)>0$, the $k$ interior roots all belong to $\Dstar$.
\end{proof}

\Needspace{11\baselineskip}
We also record the equivalent positivity statement for a real
polynomial on $[-2,2]$.

\begin{corollary}\label{cor:real}
There is a unique polynomial $R_{m,n}\in\R[x]$ of degree $k$ such that
\begin{equation}\label{eq:real-poly}
 Q_{m,n}(\zeta)=\zeta^kR_{m,n}(\zeta+\zeta^{-1}),
 \qquad\zeta\ne0.
\end{equation}
It is strictly positive on $[-2,2]$. More precisely,
\begin{equation}\label{eq:real-bound}
 R_{m,n}(x)\ge\frac m4(2-x)\quad(-2\le x\le2),
 \qquad R_{m,n}(2)=m^3.
\end{equation}
\end{corollary}

\begin{proof}
By Lemma~\ref{lem:center}, we may write
\[
 Q_{m,n}(\zeta)=\sum_{j=0}^{2k}q_j\zeta^j,\qquad
 q_j\in\R,\quad q_j=q_{2k-j},\quad q_{2k}>0.
\]
Pairing the terms equidistant from the middle coefficient gives
\begin{align*}
 \zeta^{-k}Q_{m,n}(\zeta)
 &=q_k+\sum_{j=1}^k
   \bigl(q_{k+j}\zeta^j+q_{k-j}\zeta^{-j}\bigr)\\
 &=q_k+\sum_{j=1}^kq_{k+j}(\zeta^j+\zeta^{-j}).
\end{align*}
Define real polynomials by
\[
 C_0(x)=2,\qquad C_1(x)=x,\qquad
 C_{j+1}(x)=xC_j(x)-C_{j-1}(x)\quad(j\ge1).
\]
For $x=\zeta+\zeta^{-1}$, the identities
\[
 \begin{aligned}
 C_0(x)&=2,\qquad
 C_1(x)=\zeta+\zeta^{-1},\\
 (\zeta+\zeta^{-1})(\zeta^j+\zeta^{-j})
 &-(\zeta^{j-1}+\zeta^{-(j-1)})
 =\zeta^{j+1}+\zeta^{-(j+1)}
 \end{aligned}
\]
show by induction that
\[
 C_j(\zeta+\zeta^{-1})=\zeta^j+\zeta^{-j}
 \qquad(j\ge0,\ \zeta\ne0).
\]
Moreover, $C_1(x)=x$ and $C_2(x)=x^2-2$; the recurrence then
shows that $C_j$ is monic of degree $j$ for every $j\ge1$.
Thus the polynomial
\[
 R_{m,n}(x):=q_k+\sum_{j=1}^kq_{k+j}C_j(x)
\]
satisfies
\[
 Q_{m,n}(\zeta)=\zeta^kR_{m,n}(\zeta+\zeta^{-1}),
 \qquad
 R_{m,n}(x)=q_{2k}x^k+\text{terms of lower degree}.
\]
Since $q_{2k}>0$, its degree is exactly $k$.

If another polynomial $\widetilde R$ satisfies the same identity,
substitution of $\zeta=e^{i\theta}$ gives
\[
 \widetilde R(2\cos\theta)=R_{m,n}(2\cos\theta)
 \qquad(\theta\in\R).
\]
Hence $\widetilde R-R_{m,n}$ vanishes on the whole interval
$[-2,2]$ and is identically zero. This proves uniqueness.

Finally, for $x\in[-2,2]$, choose $\theta\in[0,\pi]$ with
$x=2\cos\theta$. The identity just proved and the definition of
$F_{m,n}$ give
\[
 R_{m,n}(x)=e^{-ik\theta}Q_{m,n}(e^{i\theta})=F_{m,n}(\theta).
\]
Recalling the conclusions of Proposition~\ref{prop:strict}, we obtain
\begin{align*}
 R_{m,n}(x)&=F_{m,n}(\theta)>0,\\
 R_{m,n}(x)&\ge m\sin^2(\theta/2)
             =\frac m2(1-\cos\theta)=\frac m4(2-x),\\
 R_{m,n}(2)&=F_{m,n}(0)=m^3.
\end{align*}
The last identity also gives strict positivity at the endpoint $x=2$,
where the displayed lower bound is zero.
\end{proof}

\begin{remark}\label{rem:kfive}
The case $m-n=5$ was singled out in
\cite[Section~4.3]{EdholmMathew2026}. When $m=n+5$ and
$\gcd(n,5)=1$, Corollary~\ref{cor:real} gives a
strictly positive polynomial $R_{n+5,n}$ of degree five on $[-2,2]$.
The conclusion holds for all four residue classes
$n\equiv1,2,3,4\pmod5$ and requires no separate coefficient
calculation for these classes.
\end{remark}

\section{Proof of Corollary~\ref{cor:rational-classification}}\label{sec:zeros}

We first relate the roots of $Q_{m,n}$ to zeros of the Bergman kernel.

\begin{proposition}\label{prop:zeros}
Let $m>n$ be relatively prime positive integers. Every root
$\zeta_0\in\Dstar$ of $Q_{m,n}$ determines a pair
$(z,w)\in H_{m/n}\times H_{m/n}$ for which $K_{m/n}(z,w)=0$.
In particular, $H_{m/n}$ is not a Lu Qi-Keng domain.
\end{proposition}

\begin{proof}
Let $Q_{m,n}(\zeta_0)=0$ with $\zeta_0\in\Dstar$. Write
$\zeta_0=\rho e^{i\phi}$, where $0<\rho<1$, and set
\begin{equation}\label{eq:zero-pair}
 z=(\sqrt\rho,\sqrt\rho),\qquad
 w=(\sqrt\rho e^{-i\phi},\sqrt\rho e^{-i\phi}).
\end{equation}
Since $m>n$ and $0<\rho<1$,
\[
 |z_1|^m=|w_1|^m=\rho^{m/2}
 <\rho^{n/2}=|z_2|^n=|w_2|^n<1.
\]
This verifies that both $z$ and $w$ lie in $H_{m/n}$. Their
coordinate products are
\[
 s=z_1\overline{w_1}=\rho e^{i\phi}=\zeta_0,
 \qquad
 t=z_2\overline{w_2}=\rho e^{i\phi}=\zeta_0.
\]
Recall that $Q_{m,n}(\zeta)=\zeta^{1-2n}P_{m,n}(\zeta,\zeta)$
for $\zeta\ne0$. Hence
\[
 P_{m,n}(\zeta_0,\zeta_0)
 =\zeta_0^{2n-1}Q_{m,n}(\zeta_0)=0.
\]
The denominator of the kernel formula is nonzero, since
\begin{align*}
 |1-\zeta_0|&\ge1-\rho>0,\\
 |\zeta_0^n-\zeta_0^m|
 &=\rho^n|1-\zeta_0^{m-n}|
 \ge\rho^n(1-\rho^{m-n})>0.
\end{align*}
Evaluating the kernel formula \eqref{eq:kernel} at $s=t=\zeta_0$
therefore gives
\[
 K_{m/n}(z,w)
 =\frac{P_{m,n}(\zeta_0,\zeta_0)}
 {m\pi^2(1-\zeta_0)^2(\zeta_0^n-\zeta_0^m)^2}=0.
\]
The Edholm--Mathew conjecture, proved above, supplies exactly
$m-n\ge1$ roots in $\Dstar$, counted with multiplicity. In
particular, at least one such pair $(z,w)$ exists, and $H_{m/n}$
is not a Lu Qi-Keng domain.
\end{proof}

For the remaining rational exponents, we use the Bergman kernel
transformation law. If $\Phi:U\to V$ is a biholomorphism and
$J_\Phi=\det\Phi'$, then
\begin{equation}\label{eq:transform}
 K_U(z,w)=J_\Phi(z)
 K_V(\Phi(z),\Phi(w))\overline{J_\Phi(w)}.
\end{equation}
The chain rule gives
\[
 J_{\Phi^{-1}}(\Phi(z))J_\Phi(z)=1,
\]
so both Jacobian factors in \eqref{eq:transform} are nonzero.
Since $\Phi$ is onto, $U$ is a Lu Qi-Keng domain if and only if
$V$ is. The monomial maps used below are the ones appearing in
\cite{EdholmMcNeal2017} and
\cite[Theorem~2.19]{EdholmMathew2026}.

\begin{proof}[Proof of Corollary~\ref{cor:rational-classification}]
Suppose first that $m=1$. Define
\[
 \Psi_n:H_{1/n}\longrightarrow\D\times\Dstar,
 \qquad \Psi_n(z_1,z_2)=(z_1z_2^{-n},z_2).
\]
For $z\in H_{1/n}$,
\[
 0<|z_2|<1,\qquad |z_1z_2^{-n}|=\frac{|z_1|}{|z_2|^n}<1.
\]
Conversely, if $(u,v)\in\D\times\Dstar$, then
\[
 |uv^n|=|u||v|^n<|v|^n<1,
\]
so $(uv^n,v)\in H_{1/n}$. The maps are holomorphic, and
\[
 \Psi_n^{-1}(u,v)=(uv^n,v),\qquad
 J_{\Psi_n}(z)
 =\det\begin{pmatrix}
 z_2^{-n}&-nz_1z_2^{-n-1}\\
 0&1
 \end{pmatrix}
 =z_2^{-n}\ne0.
\]
Thus $\Psi_n$ is a biholomorphism.

For completeness, square-integrability removes the puncture from
the Bergman space of $\Dstar$. Indeed, let
\[
 f(v)=\sum_{a\in\Z}c_av^a\in A^2(\Dstar).
\]
For every integer $a$ and $0<\varepsilon<1/2$, angular orthogonality
gives
\[
 \|f\|_{L^2(\Dstar)}^2
 \ge\int_{\varepsilon}^{1/2}\int_0^{2\pi}
          |f(re^{i\vartheta})|^2r\,d\vartheta\,dr
 \ge2\pi|c_a|^2\int_{\varepsilon}^{1/2}r^{2a+1}\,dr.
\]
If $a\le-1$, the last integral tends to infinity as
$\varepsilon\downarrow0$, so $c_a=0$. Hence $f$ extends
holomorphically across $0$, without changing its $L^2$ norm.
The restriction map identifies $A^2(\D)$ and $A^2(\Dstar)$
isometrically. The product of the standard disc bases then gives
\begin{align*}
 K_{\D\times\Dstar}((u,v),(\xi,\eta))
 &=\frac1{\pi^2}\sum_{a,b\ge0}(a+1)(b+1)
              (u\overline\xi)^a(v\overline\eta)^b\\
 &=\frac1{\pi^2(1-u\overline\xi)^2(1-v\overline\eta)^2}.
\end{align*}

For $z,w\in H_{1/n}$, put
$s=z_1\overline{w_1}$ and $t=z_2\overline{w_2}$. Then
\[
 J_{\Psi_n}(z)\overline{J_{\Psi_n}(w)}=t^{-n},\qquad
 (z_1z_2^{-n})\overline{(w_1w_2^{-n})}=s/t^n.
\]
Substituting these identities into the transformation law yields
\begin{align}
 K_{1/n}(z,w)
 &=J_{\Psi_n}(z)
 K_{\D\times\Dstar}(\Psi_n(z),\Psi_n(w))
 \overline{J_{\Psi_n}(w)}\notag\\
 &=\frac{t^{-n}}{\pi^2(1-s/t^n)^2(1-t)^2}
 =\frac{t^n}{\pi^2(1-t)^2(t^n-s)^2}.
 \label{eq:reciprocal-kernel}
\end{align}
The domain inequalities imply
\[
 0<|t|<1,\qquad
 |s|=|z_1||w_1|<|z_2|^n|w_2|^n=|t|^n.
\]
Consequently,
\[
 t^n\ne0,\qquad |1-t|\ge1-|t|>0,\qquad
 |t^n-s|\ge|t|^n-|s|>0.
\]
Both the numerator and denominator of the displayed kernel are
nonzero. Thus $H_{1/n}$ is a Lu Qi-Keng domain.

Now suppose $m\ge2$. Euclidean division gives
\begin{equation}\label{eq:division}
 n=hm+r,\qquad h\in\Z_{\ge0},\qquad 1\le r\le m-1.
\end{equation}
Indeed, the initial remainder lies in $\{0,\ldots,m-1\}$, and
\[
 \gcd(m,r)=\gcd(m,n)=1.
\]
Since $m\ge2$, the remainder cannot be zero. Define
\[
 \Phi_h:H_{m/n}\longrightarrow H_{m/r},
 \qquad \Phi_h(z_1,z_2)=(z_1z_2^{-h},z_2).
\]
For $z\in H_{m/n}$, writing $u=\Phi_h(z)$ gives
\[
 |u_1|^m=|z_1|^m|z_2|^{-hm}
 <|z_2|^{n-hm}=|u_2|^r<1.
\]
Conversely, for $u\in H_{m/r}$,
\[
 |u_1u_2^h|^m=|u_1|^m|u_2|^{hm}
 <|u_2|^{r+hm}=|u_2|^n<1.
\]
The second coordinates are nonzero, so both maps are holomorphic,
and direct composition gives
\[
 \Phi_h^{-1}(u_1,u_2)=(u_1u_2^h,u_2),\qquad
 J_{\Phi_h}(z)=z_2^{-h}\ne0.
\]
Therefore $\Phi_h$ is a biholomorphism.

Because $m>r\ge1$ and $\gcd(m,r)=1$,
Proposition~\ref{prop:zeros} provides $u,v\in H_{m/r}$ with
$K_{m/r}(u,v)=0$. Set
\[
 z=\Phi_h^{-1}(u),\qquad w=\Phi_h^{-1}(v).
\]
Then $z,w\in H_{m/n}$ and, by the transformation law,
\begin{align*}
 K_{m/n}(z,w)
 &=J_{\Phi_h}(z)K_{m/r}(\Phi_h(z),\Phi_h(w))
                  \overline{J_{\Phi_h}(w)}\\
 &=z_2^{-h}K_{m/r}(u,v)\overline{w_2^{-h}}=0.
\end{align*}
Thus $H_{m/n}$ is not a Lu Qi-Keng domain when $m\ge2$.
Together with the case $m=1$, this proves the classification.
\end{proof}
\subsection*{Funding}
The project is supported by the National Natural Science Foundation of China (Grant no. 12571085), National Key Research and Development Project of China (Grant no. 2021YFA1002600) and Natural Science Foundation of Beijing Municipal (No.1252005).
\section*{Data availability}
The results are established analytically and do not rely on numerical
data.

\end{document}